\documentclass[12pt,leqno]{article}
\usepackage{amsmath,amssymb,latexsym}
\usepackage{amscd}
\usepackage{mathrsfs}
\usepackage{srcltx}
\usepackage{theorem}
\usepackage{amsfonts}

\newcommand{\R}{\mathbb{R}}        % real field
\newcommand{\N}{\mathbb{N}}        % natural numbers
\newcommand{\Rn}{\mathbb{R}^d}     % d-dim real space
\newcommand{\Cn}{\mathbb{C}^d}

\newcommand{\Di}{\cD(\Rn)}
\newcommand{\Dip}{\cD\,'(\Rn)}

\newcommand{\Proof}{\textbf{Proof:} \ }
\newcommand{\qed}{\hspace*{\fill} $\Box $}
\newcommand{\To}{\longrightarrow}

\newcommand{\cD}{\mathscr{D}}
\newcommand{\cF}{\mathscr{F}}
\newcommand{\cE}{\mathscr{E}}
\newcommand{\cS}{\mathscr{S}}
\newcommand{\F}{Fr\'echet }

\newtheorem{proposition}{Proposition}[section]

\newtheorem{lemma}[proposition]{Lemma}
\newtheorem{corollary}[proposition]{Corollary}
\newtheorem{theorem}[proposition]{Theorem}

\newcommand{\Exp}{{\rm Exp}}

\newcommand{\id}{{\rm id}}
\newcommand{\im}{{\rm im}}

\newcommand{\supp}{\mathrm{supp}\, }

\newcommand{\vp}{\varphi}

\title{{\sc Surjectivity of Euler operators on temperate distributions
}}
\author{Dietmar Vogt}
\date{}

\begin{document}

\maketitle

\footnotetext{
    \hfil\break\indent \begin{minipage}[t]{14cm}{\em Key words and phrases:} Euler differential operators, temperate distributions, global solvability, $C^\infty$-functions of exponential decay. \end{minipage}.
     \hfil\break\indent
%\centerline{\em Dedicated to }
{}}

\begin{abstract}
 Euler operators are partial differential operators of the form $P(\theta)$
where $P$ is a polynomial and $\theta_j = x_j \partial/\partial x_j$. We show
that every non-trivial Euler operator is surjective on the space of temperate
distributions on $R^d$. This is in sharp contrast to the behaviour of such
operators when acting on spaces of differentiable or analytic functions.

\end{abstract}

%\vspace{1cm}

In the present note we study Euler differential operators on the space $\cS'(\Rn)$ of temperate distributions on $\Rn$. These are operators of the form $P(\theta)$ where $P$ is a polynomial and $\theta_j = x_j\, \partial/\partial x_j$. We show that every Euler operator is surjective on $\cS'(\Rn)$ which is in sharp contrast to the behaviour in spaces of differentiable functions since the operator $P(\theta)$ is, in general, singular at the coordinate hyperplanes. Even for $d=1$ the simple example of $\theta$ acting on $C^\infty(\R)$ shows that surjectivity there is in general impossible. There are natural necessary conditions for a function to be in the range of an operator $P(\theta)$, solvability under these conditions has been shown in Doma\'nski-Langenbruch \cite{DL1}. For real analytic functions the situation is even more complicated, see \cite{DL2}. As an example our result implies the following: if $g$ is a polynomial function on $\Rn$ then the equation $P(\theta)f=g$ may not have a $C^\infty$-solution $f$ on $\Rn$
 but it will always have a temperate distribution $f$ as solution on $\Rn$.  We first study partial differential operators $P(\partial)$ with constant coefficients on the space $Y(\Rn)$ of $C^\infty$-functions with exponential decay on $\Rn$ and on its dual the space $Y(\Rn)'$ the space of distributions with exponential growth.  We show that every non-trivial operator $P(\partial)$ is surjective on $Y(\Rn)'$. By the exponential diffeomorphism this implies the surjectivity of $P(\theta)$ on the space $\cS'(Q)$ of temperate distributions on the positive quadrant on $\Rn$, hence surjectivity on $\cS'(\Rn)$ up to a distribution with support in the union of coordinate hyperplanes. By a method similar to the one used in \cite{DL1} we then show the result by induction on the dimension.

\section{Preliminaries}
We use the following notation $\partial_j=\partial/\partial x_j$, $\theta_j = x_j\,\partial_j$ and $D_j=-i\partial_j$. For a multiindex $\alpha\in\N_0^d$ we set $\partial^\alpha = \partial_1^{\alpha_1}..\partial_d^{\alpha_d}$, likewise for $\theta^\alpha$ and $D^\alpha$. For a polynomial $P(z)=\sum_\alpha c_\alpha z^\alpha$ we consider the  \it Euler operator \rm $P(\theta)=\sum_\alpha c_\alpha \theta^\alpha$ and also the operators $P(\partial)$ and $P(D)$, defined likewise.

$P(\theta)$ and $P(\partial)$ are connected in the following way.  We set for $x\in\Rn$
$$\Exp(x)=(\exp(x_1),..,\exp(x_d)).$$
$\Exp$ is a diffeomorphism from $\Rn$ onto $Q:=(0,+\infty)^d$. Therefore
$$C_\Exp:f\To f\circ\Exp$$
is a linear topological isomorphism from $C^\infty(Q)$ onto $C^\infty(\Rn)$. For $f\in C^\infty(Q)$ we have
$P(\partial)(f\circ\Exp)=(P(\theta)f)\circ\Exp$ that is $P(\partial)\circ C_\Exp =C_\Exp\circ P(\theta)$. In this way solvability properties of $P(\theta)$ on $C^\infty(Q)$ can be reduced to solvability properties of $P(\partial)$ on $C^\infty(\Rn)$. This has been done in $\cite{V1}$. We apply the same argument to the space $\cS(Q)$ where $\cS(Q)=\{f\in\cS(\Rn)\,:\,\supp f\subset \overline{Q}\}$ and $\cS(\Rn)$ is the Schwartz space of rapidly decreasing $C^\infty$-functions on $\Rn$.

Throughout the paper we use standard notation of Functional Analysis, in  particular, of distribution theory, and of the theory of partial differential operators. For unexplained notation we refer to \cite{DK}, \cite{H1}, \cite{MV}, \cite{LSI}, \cite{LS}.

\section{Distributions with exponential growth}

We start with studying partial differential operators on $\Rn$ and we will transfer our results by the exponential diffeomorphism to results on Euler operators on $Q$. We set
\begin{eqnarray*} Y(\Rn)&:=&\{f\in C^\infty(\Rn)\,:\,\sup_x|f^{(\alpha)}(x)|\,e^{k|x|}<\infty \text{ for all  }\alpha\text{ and }k\in\N\}\\
&=& \{f\in C^\infty(\Rn)\,:\,\sup_x|f^{(\alpha)}(x)|\,e^{x\eta}<\infty \text{ for all  }\alpha\text{ and }\eta\in\Rn\}
\end{eqnarray*}
with its natural topology. Here $x\eta=\sum_j x_j \eta_j$ and $|x|:=|x|_1$.

Then $Y(\Rn)$ is a \F space, closed under convolution and $P(\partial)$  is a continuous linear operator in $Y(\Rn)$ for every polynomial $P$. $\Di\subset Y(\Rn)$ as a dense subspace, hence $Y(\Rn)'\subset \Dip$. We obtain

\begin{lemma} \label{l1} $C_\Exp(\cS(Q))= Y(\Rn)$.
\end{lemma}

\Proof We first claim that
$$(f\circ \Exp)^{(\alpha)} = \sum_{\beta\le\alpha} a_\beta\, (f^{(\beta)}\circ \Exp)\, \Exp^{\beta}$$
with $a_\alpha=1$ and this is shown by induction.

This implies that for $f\in\cS(Q)$ we have
$$\sup_{x\in\Rn}|(f\circ \Exp)^{(\alpha)}(x)|\,e^{k|x|}\le \sum_{\beta\le\alpha} a_\beta \sup_{\xi\in Q} |f^{(\beta)}(\xi)| |\xi|^{|\beta|+k} <+\infty$$
for all $\alpha$ and $k\in\N$.

On the other hand we have for $g=f\circ \Exp\in Y(\Rn)$
$$(f^{(\alpha)}\circ \Exp)\, \Exp^{\alpha}=(f\circ \Exp)^{(\alpha)} - \sum_{\beta\le\alpha,\,\beta\neq\alpha} a_\beta\, (f^{(\beta)}\circ \Exp)\, \Exp^{\beta}$$
hence
$$\sup_{\xi\in Q} |f^{(\alpha)}(\xi)|\,\xi^{\alpha+\gamma}\le \sup_{x\in\Rn} |g^{(\alpha)}(x)|\, \Exp^\gamma(x) + \sum_{\beta\le\alpha,\,\beta\neq\alpha} a_\beta\, |f^{(\beta)}(\xi)|\, \xi^{\beta+\gamma}.$$
for all $\alpha\in\N_0^d$ and $\gamma\in\Rn$.
From here one derives easily by induction that $f\in\cS(\Rn)$. \qed

The space $Y(\Rn)$ can also be described by means of the Fourier transformation. This description might also be used to show Lemma \ref{t1}. We will use another method. However the description exhibits that from the point of view of the Fourier transformation $Y(\Rn)$ is a very natural space. We define
$$H_Y(\Rn)=\{g\in H(\Cn)\,:\,\sup_{x,\,|y|\le k} |x+iy|^k |g(x+iy)|<\infty \text{ for all } k\in\N\}$$
and remark that $H_Y(\Rn)\subset\cS(\Rn)$, due to Cauchy's estimates. We obtain:

\begin{proposition}\label{l2} The Fourier transformation maps $Y(\Rn)$ isomorphically onto $H_Y(\Rn)$.
\end{proposition}

\Proof For $f\in Y(\Rn)$ we have with $z=x+iy$
$$ \hat{f}(z)\,z^\alpha = (2\pi)^{-n/2} \int z^\alpha f(\xi) e^{-iz\xi}d\xi = (2\pi)^{-n/2}\,i^{|\alpha|} \int f^{(\alpha)}(\xi)\,e^{y\xi}\,e^{-ix\xi}d\xi$$
and therefore $\hat{f}\in H_Y(\Rn)$.

On the other hand, for $g\in H_Y(\Rn)$ there is $f\in\cS(\Rn)$ with $g=\hat{f}$. Then
\begin{eqnarray*}f^{(\alpha)}(x)&=&(2\pi)^{-n/2}i^{|\alpha|} \int g(\xi)\,\xi^\alpha e^{ix\xi} d\xi\\ &=& (2\pi)^{-n/2} i^{|\alpha|} \int g(\xi+i\eta)(\xi+i\eta)^\alpha e^{-x\eta} e^{ix\xi}d\xi.
\end{eqnarray*}
Therefore we have for every $\eta\in\Rn$
$$|f^{(\alpha)}(x)|\, e^{x\eta}\le (2\pi)^{-n/2}  \int |g(\xi+i\eta)(\xi+i\eta)^\alpha| d\xi <+\infty. $$
This means that $f\in Y(\Rn)$. \qed

Our first main result is based on the following Lemma.
\begin{lemma}\label{t1} For every non-trivial polynomial the operator $P(-\partial)$ in $Y(\Rn)$ is an isomorphism onto its range.
\end{lemma}

\Proof  We use the construction of an elementary solution of P. Wagner \cite[Proposition 1]{W} and assume $P$ to be nontrivial of degree $m$. For $\eta\in\Rn$ with $P_m(\eta)\neq0$ and pairwise different real numbers $\lambda_0,..,\lambda_m$ we set with suitable coefficients $a_j$
$$E=\frac{1}{P_m(2\eta)}\sum_{j=0}^m a_j e^{\lambda_j\eta x} \cF_\xi^{-1}\left(\frac{\overline{P(i\xi + \lambda_j \eta}}{P(i\xi + \lambda_j \eta}\right).$$
Then, by Wagner, loc. cit., $E$ is an elementary solution for $P(\partial)$, that is, $\langle E, P(-\partial)\vp)\rangle = \vp(0)$ for any $\vp\in \Di$. By continuous extension this holds also for $\vp\in Y(\Rn)$ since $E\in Y(\Rn)'$. That means that for any $\vp\in Y(\Rn)$ the term
$\vp(y)$ is a linear combination of terms
$$\Psi_j(y)=\left\langle (P(-\partial_x)\vp(x+y))e^{\lambda_j \eta x}, \cF_\xi^{-1}\left(\frac{\overline{P(i\xi + \lambda_j \eta}}{P(i\xi + \lambda_j \eta}\right)\right\rangle.$$
Then $\Psi_j(y) e^{\lambda_j\eta y}$ has the form
$$\Psi_j(y) e^{\lambda_j\eta y}=\left\langle \psi(\cdot+y), \cF^{-1}(G)\right\rangle$$
where $\psi(x)= (P(-\partial_x)\vp)(x))e^{\lambda_j \eta x} \in Y(\Rn)$ and $G\in \L_\infty(\Rn)$, $\|G\|_\infty =1$. We obtain
$$|\Psi_j(y) e^{\lambda_j\eta y}|=\left|\int \hat{\psi}(-\xi) e^{i\xi y} G(\xi) d\xi\right|\le \int |\hat{\psi}(\xi)|\, d\xi= \|P(-\partial)\vp\|_{Y(\Rn)}$$
where $\|\cdot\|_{Y(\Rn)}$ is a semi-norm in $Y(\Rn)$ (cf. Proposition \ref{l2}). We do that for every $j$ and find a semi-norm $\|\cdot\|_\eta$ in $Y(\Rn)$ such that
$$|\vp(y)|\le \frac{1}{|P_m(2\eta)|}\sum_{j=0}^m |a_j|\,e^{-\lambda_j \eta y}\|P(-\partial)\vp\|_\eta.$$
We may assume that our choice was so that $\lambda_j\ge 1$ for all $j$. We evaluate the inequality separately for every quadrant
$Q_e=\{x\,:\,e_j x_j\ge 0\text{ for all }j\}$, $e_j=\pm 1$ for all $j$. For given $c>0$ we  choose $\eta_e\in Q_e$ such that $\eta_e y\ge c\,|y|$ for all $y\in Q_e$.
This yields with a proper constant $D$
$$|\vp(y)|\le D \max_e \|P(-\partial)\vp\|_{\eta_e} e^{-c|y|}=: \|P(-\partial)\vp\|\, e^{-c|y|}.$$
We may apply this estimate to $\vp^{(\alpha)}$ for any $\alpha$. Since $\partial^\alpha$ commutes with $P(-\partial)$ and since $c>0$ was chosen arbitrarily,  we obtain the result. \qed

%$$\left\langle \cF \Big(P(-\partial_x)(\vp(x+y))e^{\lambda_j \eta x}\Big), \frac{\overline{P(i\xi + \lambda_j \eta}}{P(i\xi + \lambda_j %\eta}\right\rangle.$$

As an immediate consequence we obtain by dualization of Lemma \ref{t1}.

\begin{theorem}\label{t3} Every non-trivial operator $P(\partial)$ is surjective on $Y(\Rn)'$.
\end{theorem}

\section{Euler differential operators on $\cS'(\Rn)$}

From Lemma \ref{t1} we derive the following lemma which is the basis for our main result.

\begin{lemma} \label{c1} Every non-trivial Euler operator is surjective on $\cS'(Q)$.
\end{lemma}

\Proof To show surjectivity of $P(\theta)$ on $\cS'(Q)$ we have to show that $P(\theta^*)$ is an isomorphism onto its range in $\cS(Q)$. Here $\theta^*$ is the transpose of $\theta$, that is, $\theta_j^*\vp=-\vp-\theta_j\vp$ and therefore $P(\theta^*)=P(-1-\theta)$.
By Lemma \ref{l1} and the fact that  $P(-1-\partial)\circ C_\Exp =C_\Exp\circ P(\theta^*)$, we obtain from Lemma \ref{t1} that $P(\theta^*)$ is an isomorphism onto its range in $\cS(Q)$ which shows the result. \qed

We set $Z_0=\{x\in\Rn\,:\,\exists j : x_j=0\}$ and $\cS_0=\{f\in\cS(\Rn)\,:\, f\text{ is flat on }Z_0\}$.
Since Lemma \ref{c1} holds, mutatis mutandis, on every `quadrant' we obtain:

\begin{corollary}\label{c2} Every non-trivial Euler operator is surjective on $\cS_0'$.
\end{corollary}

Let $T\in\cS'(\Rn)$ and $P(\theta)$ be a non-trivial Euler operator. We want to solve the equation $P(\theta)U=T$. We set $T_0=T|_{\cS_0}$. Then $T_0\in\cS_0'$ and we can find $U_0\in\cS_0'$ such that $P(\theta)U_0=T_0$. We extend $U_0$ by use of the Hahn-Banach Theorem to a distribution $U_1\in\cS'(\Rn)$ and for any solution $U$ of our problem we have $P(\theta)(U-U_1)=0$ on $\cS_0$. Hence it suffices to solve the equation $P(\theta)U_2=T_1$ with $T_1=T-P(\partial)U_1\in \cS'(Z_0)=\{S\in\cS'(\Rn)\,:\; S\text{ vanishes on }\cS_0\}$.

By $H_j$ we denote the coordinate hyperplanes $H_j=\{x\in\Rn\,:\,x_j=0\}$. We set $\cS_{j_1,..,j_k}=\{f\in\cS(\Rn)\,:\,f\text{ flat on } H_{j_1}\cup..\cup H_{j_k}\}$, then $\cS_0=\cS_{1,..,d}$, and we set $\cS'(H_{j_1}\cup .. \cup H_{j_k})=\{S\in\cS'(\Rn)\,:\; S\text{ vanishes on }\cS_{j_1,..,j_k}\}$.

\begin{lemma}\label{l3} Every $T\in\cS'(Z_0)$ has a decomposition $T=T_1+\dots+T_d$ where $T_j\in\cS'(H_j)$ for all $j$.
\end{lemma}

\Proof We act by induction and consider the canonical map
$$\Phi_k:=\cS(\Rn)/\cS_{k,..,d}\to \cS(\Rn)/\cS_k \times \cS(\Rn)/\cS_{k+1,..,d}.$$
It is injective and we want to show that its image is closed. We admit that $\cS_k+\cS_{k+1,..,d}$ is closed in $\cS(\Rn)$ and we will show this later. We consider the map
$$\Psi_k: \cS(\Rn)/\cS_k \times \cS(\Rn)/\cS_{k+1,..,d} \to \cS(\Rn)/\cS_k+\cS_{k+1,..,d}$$
given by
$$(\vp_1 + \cS_k)\times (\vp_2 + \cS_{k+1,..,d})\mapsto (\vp_1-\vp_2)+(\cS_k+\cS_{k+1,..,d}).$$
We claim that $\im \Phi_k=\ker \Psi_k$. One inclusion is evident, for the other assume that $(\vp_1-\vp_2) \in (\cS_k+\cS_{k+1,..,d})$. Then there are $\psi_1\in\cS_k$ and $\psi_2\in\cS_{k+1,..,d}$ such that $\vp_1-\vp_2=\psi_1-\psi_2$, hence $\vp_1-\psi_1=\vp_2-\psi_2=:\vp$. So we have $\Phi_k(\hat{\vp})= \hat{\vp_1}\times\Hat{\vp_2}$, where $\hat{\ }$ denotes the respective equivalence class.

Let $T_{j_1,..,j_k}$ denote a distribution in $\cS'(H_{j_1}\cup .. \cup H_{j_k})$. Dualization of the previous shows that for every $T_{k,..,d}$ there  are $T_k$ and $T_{k+1,..,d}$ such that $T_{k,..,d}=T_k + T_{k+1,..,d}$. Starting with $T=T_{1,..,d} \in\cS'(Z_0)$, induction over $k$ yields the result.

It remains to show that $\cS_k+\cS_{k+1,..,d}$ is closed in $\cS(\Rn)$. We write $\cS(\Rn)=\cS(\R^k)\widehat{\otimes}_\pi\cS(\R^{d-k})$ and consider the map $$A_1\otimes A_2: \cS(\R^k)\widehat{\otimes}_\pi\cS(\R^{d-k})\to \cS(\R^k)/S_k(\R^k)\widehat{\otimes}_\pi\cS(\R^{d-k})/S_{k+1,..,d}(\R^{d-k}).$$
Here $A_1,\,A_2$ are the respective quotient maps and $S_k(\R^k)$ resp. $S_{k+1,..,d}(\R^{d-k})$ are self-explaining. By use of Lemma \ref{l4} below we get
$$\ker (A_1\otimes A_2) =  S_k(\R^k)\widehat{\otimes}_\pi \cS(\R^{d-k}) +\cS(\R^k)\widehat{\otimes}_\pi S_{k+1,..,d}(\R^{d-k})
=\cS_k+\cS_{k+1,..,d}$$
and the claim is shown.  \qed

The following lemma is known (see \cite{DL2}, Theorem 2.12). We sketch a proof for the convenience of the reader.

\begin{lemma}\label{l4} Let $A_1:E_1\to F_1$,\,$A_2:E_2\to F_2$ be continuous linear maps between nuclear \F spaces then for
$A_1\otimes A_2 : E_1\widehat{\otimes}_\pi E_2\to F_1\widehat{\otimes}_\pi F_2$ we have
$\ker(A_1\otimes A_2)=\ker A_1\widehat{\otimes}_\pi E_2 + E_1\widehat{\otimes}_\pi \ker A_2$.
\end{lemma}

\Proof  By Grothendieck \cite{G}, Chap I,\,p.\,38 and Chap II,\,p.\,70, we obtain $\ker\, (A_1\otimes \id_{E_2})=\ker A_1\widehat{\otimes}_\pi E_2$ and $\ker\, (\id_{F_1}\otimes A_2)=F_1\widehat{\otimes}_\pi \ker A_2$. We calculate the kernel of the composition, which is $A_1\otimes A_2$, and obtain the result. \qed

%\vspace{1cm}

We are now ready to prove our main result. We will use the structure of distributions in $\cS'(H_j)$ (see in analogy \cite[Chap III, Th\'eor\`eme XXXVI]{LSI}) and the fact that $\partial_j$ and multiplication with $x_j$, which is up to a factor the Fourier transform of $\partial_j$, are surjective on $\cS'(\Rn)$.
\begin{theorem}\label{t2} Every non-trivial Euler operator is surjective on $\cS'(\Rn)$.
\end{theorem}

\Proof
The proof will be by induction on the dimension. We start with the induction step. Assume that $d>1$ and the result is shown for $d-1$.
We may assume that $P$ is irreducible.

Due to Lemma \ref{l3} and the argument before this Lemma it is enough to assume that $T\in\cS'(H_j)$. We may assume that $j=1$. We set $x'=(x_2,\dots,x_d)$. Then
$$T(\vp)= \sum_k \delta^{(k)}(x_1) T_k(x') \vp(x_1,x').$$
and analogously for $U$ if we try to solve $P(\theta)U=T$ with $U\in \cS'(H_1)$.  Here $U_k,\, T_k\in \cS'(\R^{d-1})$ and the sums are finite.

We recall that for $\psi\in C^\infty(\R)$ we have
$$\langle \theta(\delta^{(k)}),\psi\rangle=(-1)^{k+1}(x\psi)^{(k+1)}(0) = (-1)^{k+1} (k+1)\psi^{(k)}(0)=-(k+1)\langle \delta^{(k)},\psi\rangle,$$
that is, $\theta\delta^{(k)}=(-k-1)\delta^{(k)}$.

For $P(z)=\sum_\alpha c_\alpha z^\alpha$ we obtain
$$P(\theta) U =\sum_k\sum_\alpha c_\alpha (-k-1)^{\alpha_1}\delta^{(k)}{\theta'}^{\alpha'} U_k=\sum_k \delta^{(k)}P(-k-1,\theta')U_k$$
and the equation we want to solve takes the form
$$\sum_k \delta^{(k)}P(-k-1,\theta') U_k =\sum_k \delta^{(k)} T_k.$$
By the induction assumption $P(-k-1,\theta') U_k=T_k$ is solvable if $P(-k-1,z')\not\equiv 0$.

 So we have to consider the case where $P(-k-1,z')\equiv 0$. Then $z_1+k+1$ divides $P$. Since $P$ is irreducible $P(z)=C\,(z_1+k+1)$ for some constant $C\neq0$.

 Since under the canonical isomorphism $\cS'(\Rn)=\cS'(\R)\widehat{\otimes}_\pi\cS'(\R^{d-1})$ the operator $\theta_1+n$ in $\cS'(\Rn)$ corresponds to $(\theta+n) \otimes \id$ in $ \cS'(\R)\widehat{\otimes}_\pi\cS'(\R^{d-1})$ we are, due to Grothendieck's exactness theorem, reduced to the surjectivity of $\theta +n$ in $\cS'(\R)$ which is shown below.

For $d=1$ we have to consider $P(z)=z-a$. It is enough to find, for all $a$ and $k$, a distribution $S\in\cS'(\R)$ such that $(\theta-a)S=\delta^{(k)}$. Since $(\theta - a)\,\delta^{(k)}= -(k+1+a)\,\delta^{(k)}$ this is evident for $a\neq -(k+1)$. For the equation $(\theta+k+1)S=\delta^{(k)}$ we obtain a solution as follows: choose any $U$ with $x^{k+1} U = Y(x)$, where
$Y(x)=1$ for $x\ge 0$, $Y(x)=0$ otherwise, then $x^k((k+1)U+\theta U)=\delta$, hence $(\theta+k+1)U= (-1)^k \,1/k!\, \delta^{(k)}+\sum_{j=0}^{k-1}c_j\,\delta^{(j)}$ with suitable $c_j$. This is easily checked on monomials. Therefore  $S= (-1)^k\,k! \, \big(U-\sum_{j=0}^{k-1}\frac{c_j}{k-j} \, \delta^{(j)}\big)$ is a solution.\qed

Let us finally remark that for $d=1$, that is, for the case of an ordinary Euler differential operator $P(\theta)$ the proof of surjectivity can, by the fundamental Theorem of algebra, be reduced to the case $P(\theta)=\theta-a$ and be carried out in a much more direct way. The proof then shows that $P(\theta)$ is surjective also in the space $\cD'(\R)$ of Schwartz distributions. Whether this is true also for higher dimensions is not known and it is an open problem.

%\begin{lemma}\label{l3} $Y(\Rn)$ is closed under convolution.
%\end{lemma}

%\Proof let $f,g\in Y(\Rn)$ then
%$$(f*g)^{(\alpha)}(x)e^{x\eta} = \int f^{(\alpha)}(\xi)e^{\xi\eta} g(x-\xi)e^{(x-\xi)\eta} d\xi$$
%wich implies the result since $f^{(\alpha)}(\xi)e^{\xi\eta}$ and $g(\xi)e^{\xi\eta}$ are in $L_1$ due to the assumption. \qed

%For $T\in Y(\Rn)'$ and $f\in Y(\Rn)$ we have for $C$ and $\beta$ depending on $T$
%$$|(T_\xi f(x-\xi))^{(\alpha)}|\, e^{x\eta}\le C \sup_\xi |f^{(\alpha)}(x-\xi)| e^{\xi\beta+x\eta}$$

%Let $T\in Y(\Rn)'$ and $S\in Y(\Rn)'$, $f\in Y(\Rn)$ then we have with $C,\,\alpha,\eta$ chosen for $T$, $D,\,\beta, \eta'$ chosen for $S$
%$$|(T*S)f| = |T_x(S_y f(x+y))|\le C\sup_x |S_y f^{(\alpha)}(x+y)|\,e^{x\eta}\le CD\sup_{x,\,y} |f^{(\alpha+\beta)}(x+y)|e^{x\eta}e^{y\eta'}$$

\vspace{.5cm}

\noindent Bergische Universit\"{a}t Wuppertal,
\newline Dept. of Math., Gau\ss -Str. 20,
\newline D-42119 Wuppertal, Germany
\newline e-mail: dvogt@math.uni-wuppertal.de

\end{document}